\documentclass[12pt]{article}
\usepackage{amsmath,amsfonts,amssymb,amsthm}
\begin{document}

\newcommand{\1}{{{\bf 1}}}
\newcommand{\id}{{\rm id}}
\newcommand{\Hom}{{\rm Hom}\,}
\newcommand{\End}{{\rm End}\,}
\newcommand{\Res}{{\rm Res}\,}
\newcommand{\Image}{{\rm Im}\,}
\newcommand{\Ind}{{\rm Ind}\,}
\newcommand{\Aut}{{\rm Aut}\,}
\newcommand{\Ker}{{\rm Ker}\,}
\newcommand{\gr}{{\rm gr}}
\newcommand{\Der}{{\rm Der}\,}
\newcommand{\Z}{\mathbb{Z}}
\newcommand{\Q}{\mathbb{Q}}
\newcommand{\C}{\mathbb{C}}
\newcommand{\N}{\mathbb{N}}
\newcommand{\g}{\mathfrak{g}}
\newcommand{\gl}{\mathfrak{gl}}
\newcommand{\h}{\mathfrak{h}}
\newcommand{\wt}{{\rm wt}\,}
\newcommand{\A}{\mathcal{A}}
\newcommand{\D}{\mathcal{D}}
\newcommand{\Lie}{\mathcal{L}}
\newcommand{\E}{\mathcal{E}}

\def \b{\beta}
\def \<{\langle}
\def \>{\rangle}
\def \be{\begin{equation}\label}
\def \ee{\end{equation}}
\def \bex{\begin{exa}\label}
\def \eex{\end{exa}}
\def \bl{\begin{lem}\label}
\def \el{\end{lem}}
\def \bt{\begin{thm}\label}
\def \et{\end{thm}}
\def \bp{\begin{prop}\label}
\def \ep{\end{prop}}
\def \br{\begin{rem}\label}
\def \er{\end{rem}}
\def \bc{\begin{coro}\label}
\def \ec{\end{coro}}
\def \bd{\begin{de}\label}
\def \ed{\end{de}}

\newtheorem{thm}{Theorem}[section]
\newtheorem{prop}[thm]{Proposition}
\newtheorem{coro}[thm]{Corollary}
\newtheorem{conj}[thm]{Conjecture}
\newtheorem{exa}[thm]{Example}
\newtheorem{lem}[thm]{Lemma}
\newtheorem{rem}[thm]{Remark}
\newtheorem{de}[thm]{Definition}
\newtheorem{hy}[thm]{Hypothesis}
\makeatletter \@addtoreset{equation}{section}
\def\theequation{\thesection.\arabic{equation}}
\makeatother \makeatletter

\begin{Large}
\begin{center}
\textbf{Classification of irreducible modules of the vertex algebra $V_L^+$ when $L$ is a nondegenerate even lattice of an arbitrary rank}
\end{center}
\end{Large}
\begin{center}{
Gaywalee Yamskulna\footnote{E-mail: gyamsku@ilstu.edu, Partially supported by a Pre-Tenure Faculty Initiative grant from the College of Arts and Sciences, ISU}\\
Department of Mathematical Sciences, Illinois State University, Normal, IL 61790\\
and\\
Institute of Science, Walailak University, Nakon Si Thammarat,
Thailand}
\end{center}
\begin{abstract}
In this paper, we first classify all irreducible modules of the vertex algebra $V_L^+$ when $L$ is a negative definite even lattice of arbitrary rank. In particular, we show that any irreducible $V_L^+$-module is isomorphic to a submodule of an irreducible twisted $V_L$-module. We then extend this result to a vertex algebra $V_L^+$ when $L$ is a nondegenerate even lattice of finite rank.
\end{abstract}
Key Words: Vertex algebra.
\section{Introduction}

It is well known that for a nondegenerate even lattice $L$, there is a corresponding vertex algebra, $V_L$, with an automorphism $\theta$ of order 2 which is induced from the -1 isometry of the lattice (cf. \cite{b, flm}). Moreover, the $\theta$-invariant vertex sub-algebra $V_L^+$ is irreducible (cf. \cite{dm}).  
These vertex algebras $V_L^+$ provide a large class of concrete and important examples of vertex algebras. In fact, the study of $V_L^+$ was initiated in \cite{flm} during the course of the construction of the moonshine module.

To a large extent, the study of vertex algebras is the study of their representations. As for classical algebras, the study of complete reducibility of modules is of fundamental importance. The classification of irreducible $V_L^+$-modules and the complete reducibility of $V_L^+$-modules for the case when $L$ is a positive definite even lattice of an arbitrary rank and when $L$ is a rank one negative definite even lattice has been done by Abe, Dong, Jiang, Jordan and Nagatomo (cf. \cite{a,ad,dj, dn1, j}).

In this paper, we take a further step in understanding $V_L^+$ by studying classifications of irreducible $V_L^+$-modules when $L$ is a negative definite even lattice of a finite rank and when $L$ is a nondegenerate  even lattice of a finite rank. For these two cases, we prove that any irreducible $V_L^+$-module is isomorphic to a submodule of an irreducible twisted $V_L$-module. The main idea of the proofs is to classify irreducible modules of the Zhu algebras $A(V_L^+)$. There is a one-to-one correspondence between the set of inequivalent irreducible $A(V_L^+)$-modules and the set of inequivalent irreducible (admissible) $V_L^+$-modules (cf. \cite{z}).

This paper is organized as follows. In Section 2, we review the definition of a vertex algebra, and recall various notions of its twisted modules. In addition, we discuss a Zhu algebra and its properties. Also, we recall the constructions of $M(1)^+$, $V_L^+$, their modules and related topics that we will need in later sections. In Section 3, we study the Zhu algebra of $V_L^+$ when $L$ is a negative definite even lattice of a finite rank and we use this information to classify all irreducible (admissible) $V_L^+$-modules. In Section 4, we classify all irreducible (admissible) $V_L^+$-modules when $L$ is a nondegenerate even lattice of a finite rank.

{\bf Acknowledgments:} We want to thank Chongying Dong for reading this manuscript. In this research, we greatly benefited from studying the work of Abe and Dong in \cite{ad}. We have adopted many of their important ideas into this paper.
\section{Preliminaries}
\subsection{Vertex algebras and Zhu algebras}

First, we define vertex algebras, and their automorphisms. Next, we recall various notions of twisted modules for a vertex algebra. We also discuss about Zhu algebras.

\begin{de}\cite{lli} A {\em vertex algebra} $V$ is a vector space equipped with a linear map
$Y(\cdot, z):V\rightarrow(\End V)[[z,z^{-1}]], v\mapsto Y(v,z)=\sum_{n\in\Z}v_nz^{-n-1}$ and a distinguished vector ${\bf 1}\in V$ which satisfies the following properties:
\begin{enumerate}
\item $u_nv=0$ for $n>>0$.
\item $Y({\bf 1},z)=id_V$.
\item $Y(v,z){\bf 1}\in V[[z]]$ and $\lim_{z\rightarrow 0}Y(v,z){\bf 1}=v$.
\item (the Jacobi identity) \begin{eqnarray*}
& &z_0^{-1}\delta\left(\frac{z_1-z_2}{z_0}\right)Y(u,z_1)Y(v,z_2)-z_0^{-1}\delta\left(\frac{z_2-z_1}{-z_0}\right)Y(v,z_2)Y(u,z_1)\\
&&=z_2^{-1}\delta\left(\frac{z_1-z_0}{z_2}\right)Y(Y(u,z_0)v,z_2).
\end{eqnarray*}
\end{enumerate}
We denote the vertex algebra just defined by $(V, Y, {\bf 1})$ or, briefly, by $V$.
\end{de}
\begin{de} A $\Z$-graded vertex algebra is a vertex algebra $$V=\oplus_{n\in\Z}V_n; \text{ for }v\in V_n,\ \  n=\wt v,$$ equipped with a conformal vector $\omega\in V_2$ which satisfies the following relations:
\begin{itemize}
\item $[L(m),L(n)]=(m-n)L(m+n)+\frac{1}{12}(m^3-m)\delta_{m+n,0}c_V$ for $m,n\in \Z$, where $c_V\in \C$ (the central charge) and $$Y(\omega,z)=\sum_{n\in\Z}L(n)z^{-n-2}\left(=\sum_{m\in\Z}\omega_mz^{-m-1}\right);$$
\item $L(0)v=nv=(\wt v) v$ for $n\in\Z$, and $v\in V_n$;
\item $Y(L(-1)v,z)=\frac{d}{dz}Y(v,z)$.
\end{itemize}
\end{de}
\begin{de}\cite{lli} A vertex sub-algebra of a vertex algebra $V$ is a vector sub-space $U$ of $V$  such that ${\bf 1}\in U$ and  such that $U$ is itself a vertex algebra.
\end{de}

\begin{de}\cite{lli} An {\em automorphism} of a vertex algebra $V$ is a linear isomorphism of $V$ such that $g({\bf 1})={\bf 1}$, and $gY(v,z)u=Y(g(v),z)g(u)$ for $u,v\in V$.
\end{de}

Let $G$ be a finite automorphism group of $V$. We denote the $G$-fixed-point vertex sub-algebra of $V$ by $V^G(=\{v\in V | gv=v \text{ for all } g\in G\})$. 

For $g$, an automorphism of the vertex algebra $V$ of order $T$, we denote the decomposition of $V$ into eigenspaces with respect to the action of $g$ as $V=\oplus_{r=0}^{T-1}V^r$, where $V^r=\{v\in V|g(v)=e^{2\pi i r/T}v\}$.
\begin{de}\cite{d,ffr,flm,l} A {\em weak $g$-twisted $V$-module} $M$ is a vector space equipped with a linear map $Y_M(\cdot, z):V\rightarrow(\End  M)\{z\}, v\mapsto Y_M(v,z)=\sum_{n\in\Q}v_nz^{-n-1}$ which satisfies the following properties: for $v\in V$, $u\in V^r$, and $w\in M$
\begin{enumerate}
\item $v_nw=0$ for $n>>0$.
\item $Y_M(u,z)=\sum_{n\in\Z}u_{n+\frac{r}{T}}z^{-n-\frac{r}{T}-1}$.
\item $Y_M({\bf 1},z)=id_M$.
\item (the twisted Jacobi identity) \begin{eqnarray*}
& &z_0^{-1}\delta\left(\frac{z_1-z_2}{z_0}\right)Y_M(u,z_1)Y_M(v,z_2)-z_0^{-1}\delta\left(\frac{z_2-z_1}{-z_0}\right)Y_M(v,z_2)Y_M(u,z_1)\\
&&=z_2^{-1}\left(\frac{z_1-z_0}{z_2}\right)^{-r/T}\delta\left(\frac{z_1-z_0}{z_2}\right)Y_M(Y(u,z_0)v,z_2).
\end{eqnarray*}
\end{enumerate}
If $g$ is the identity map, a weak $g$-twisted $V$-module is called a weak $V$-module.
\end{de}

\begin{rem} Let $G$ be a finite automorphism group of $V$ and let $g$ be a member of $G$. Then any $g$-twisted weak $V$-module is a weak $V^G$-module.
\end{rem}

\begin{de} An {\em irreducible} weak $g$-twisted $V$-module is a weak $g$-twisted $V$-module that has no weak $g$-twisted $V$-submodule except 0 and itself. Here, a weak $g$-twisted sub-module is defined in the obvious way.
\end{de}
\begin{de}\cite{dlm} An ({\em ordinary}) $g$-twisted $V$-module is a weak $g$-twisted $V$-module $M$ which carries a $\C$-grading induced by the spectrum of $L(0)$. Then $M=\oplus_{\lambda\in\C}M_{\lambda}$ where $M_{\lambda}=\{w\in M| L(0)w=\lambda w\}$, and $\dim M_{\lambda} <\infty$. Moreover, for fixed $\lambda$, $M_{\frac{n}{T}+\lambda}=0$ for all small enough integers $n$.
\end{de}
\begin{de}\cite{dlm} An {\em admissible} $g$-twisted $V$-module $M$ is a $\frac{1}{T}\Z_{\geq 0}$-graded weak $g$-twisted $V$-module $M=\oplus_{n\in\frac{1}{T}\Z_{\geq 0}}M(n)$ such that $$v_mM(n)\subset M(n+\wt v-m-1)$$ for any homogeneous $v\in V$ and $m,n\in\Q$. Here, $\Z_{\geq 0}$ is the set of nonnegative integers.

An admissible $g$-twisted $V$-submodule of $M$ is a weak  $g$-twisted $V$-submodule $N$ of $M$ such that $N=\oplus_{n\in \frac{1}{T}\Z_{\geq 0}}N\cap M(n)$.
\end{de}
If $g$ is the identity map, then an admissible $g$-twisted $V$-module is called an admissible $V$-module.
\begin{lem}\cite{dlm} Any ordinary $g$-twisted $V$-module is an admissible $g$-twisted $V$-module.\end{lem}
\begin{de} An {\em irreducible admissible} $g$-twisted $V$-module is an admissible $g$-twisted $V$-module that has no admissible $g$-twisted submodule except 0 and itself.
\end{de}

Next, we will define a Zhu algebra and we will discuss about the relationship between the Zhu algebra and a vertex algebra. 

Let $V$ be a $\Z$-graded vertex algebra. For a homogeneous vector $u\in V$, $v\in V$, we define products $u*v$, and $u\circ v$ as follow:
\begin{eqnarray*}
u*v&=&\Res_{z}\left(\frac{(1+z)^{\wt u}}{z} Y(u,z)v\right)\\
u\circ v&=&\Res_{z}\left(\frac{(1+z)^{\wt u}}{z^2}Y(u,z)v\right).
\end{eqnarray*}  Then we extend these products linearly on $V$. We let $O(V)$ be the linear span of $u\circ v$ for all $u,v\in V$ and we set $A(V)=V/O(V)$. Also, for $v\in V$, we denote $v+O(V)$ by $[v]$.
\begin{thm}\label{zhu}\cite{z}\ \
\begin{enumerate}
\item The product $*$ induces the structure of an associative algebra on $A(V)$ with the identity ${\bf 1}+O(V)$. Moreover, $\omega+O(V)$ is a central element of $A(V)$.
\item The map $M\rightarrow M(0)$ gives a bijection between the set of equivalence classes of irreducible admissible $V$-modules to the set of equivalence classes of simple $A(V)$-modules.
\item Assume that $u\in V$ is homogeneous, $v\in V$ and $n\geq 0$. Then $$\Res_{z}\left(\frac{(1+z)^{\wt u}}{z^{2+n}}Y(u,z)v\right)\in O(V).$$
\item Let $v\in V$ be homogeneous and $u\in V$. Then $$u*v-\Res_{z}\left(\frac{(1+z)^{\wt v-1}}{z}Y(v,z)u\right)\in O(V).$$
\item For homogeneous vectors $u,v\in V$, we have $$u*v-v*u-\Res_{z}(1+z)^{\wt u -1}Y(u,z)v\in O(V).$$
\item For any $n\geq 1$, $L(-n)\equiv (-1)^n\{(n-1)((L(-2)+L(-1))+L(0)\}\mod O(V)$ where $L(n)$ are the Virasoro operators given by $Y(\omega,z)=\sum_{n\in \Z}L(n)z^{-n-2}$.
\end{enumerate}
\end{thm}
\subsection{Vertex algebras $V_L^+$ and $M(1)^+$ }

First, we will recall the construction of a vertex algebra $M(1)$ and its modules. Next, we will discuss about a vertex algebra $V_L$ where $L$ is an even lattice equipped with a nondegenerate symmetric $\Z$-bilinear form. In addition, we will review constructions of certain twisted modules of a vertex algebra $M(1)$ and a vertex algebra $V_L$. Finally, we will discuss about a vertex algebra $M(1)^+$, a vertex algebra $V_L^+$ and their irreducible modules.

Let $L$ be a rank $d$ even lattice equipped with a nondegenerate symmetric $\Z$-bilinear form $\<\cdot,\cdot\>$. We set $\h=\C\otimes_{\Z}L$ and extend $\<\cdot,\cdot\>$ to a $\C$-bilinear form on $\h$. Let $\hat{\h}=\h\otimes \C[t,t^{-1}]\oplus\C c$ be the affinization of $\h$ with the following commutator formula:
$$[\beta\otimes t^m,\alpha\otimes t^n]=m\<\beta,\alpha\>\delta_{m,-n}c\ \ \text{and}\ \ [c,\hat{h}]=0$$ for any $\alpha, \beta\in\h$, $m,n\in\Z$. Then $\hat{\h}^+=\h\otimes \C[t]\oplus\C c$ is a subalgebra of $\hat{\h}$. Let $\lambda\in \h$, and consider the  induced $\hat{\h}$-module
$$M(1,\lambda)=U(\hat{\h})\otimes_{U(\hat{\h}^+)}\C_{\lambda}$$
where $\h\otimes t\C[t]$ acts trivially on $\C$, $h$ acts as $\<h,\lambda\>$ for $h\in\h$, and $c$ acts on $\C$ as a multiplication by 1. We set $M(1)=M(1,0)$. For $\alpha\in \h$ and $n\in \Z$, we write $\alpha(n)$ for the operator $\alpha\otimes t^n$ and set $\alpha(z)=\sum_{n\in\Z}\alpha(n)z^{-n-1}.$ For $\alpha_1,...,\alpha_k\in\h$, $n_1,...,n_k>0$ and $v=\alpha_1(-n_1)...\alpha_k(-n_k)\otimes 1\in M(1)$, we define
$$Y(v,z)=: \partial^{(n_1-1)}\alpha_1(z)...\partial^{(n_k-1)}\alpha_k(z):,$$
where $\partial^{(n)}=\frac{1}{n!}\left(\frac{d}{dz}\right)^n$ and the normal ordering $:\cdot :$ is an operation which reorders the operators so that $\beta(n)$ ($\beta\in\h,n<0$) to be placed to the left of $\gamma(n)$ ($\gamma\in\h, n\geq 0$) before the expression is evaluated. We extend $Y$ to all $v\in M(1)$ by linearity. Let $\{\beta_1,...,\beta_d\}$ be an orthonormal basis of $\h$. Set ${\bf 1}=1$ and $\omega=\frac{1}{2}\sum_{i=1}^d\beta_i(-1)^2\otimes 1$. 
\begin{prop}\cite{gu} $(M(1), Y, {\bf 1},\omega)$ is a $\Z$-graded vertex algebra and $\{ M(1,\lambda) |\lambda\in\h\}$ is the set of all inequivalent irreducible $M(1)$-modules.
\end{prop}

Let $\hat{L}$ be the canonical central extension of $L$ by the cyclic group $\<\kappa \>$ of order 2. Let $e:L\rightarrow \hat{L}$ be a section such that $e_0=1$ and $\epsilon:L\times L\rightarrow\<\kappa\>$ be the corresponding 2-cocycle. We can assume that $\epsilon$ is bimultiplicative. Then $\epsilon(\alpha,\beta)/\epsilon(\beta, \alpha)=\kappa^{\<\alpha,\beta\>}$, $\epsilon(\alpha,\beta)\epsilon(\alpha+\beta, \gamma)=\epsilon(\beta,\gamma)\epsilon(\alpha, \beta+\gamma),$ and $e_{\alpha}e_{\beta}=\epsilon(\alpha,\beta)e_{\alpha+\beta}$ for $\alpha,\beta, \gamma\in L$.  We define
$$\C[L]=\oplus_{\alpha\in L}\C e^{\alpha}.$$ The action of $\hat{L}$ on $\C[L]$ is given by
$e_{\alpha}e^{\beta}=\epsilon(\alpha,\beta)e^{\alpha+\beta}$, and $\kappa e^{\beta}=-e^{\beta}$. Also, we define an action of $\h$ on $\C[L]$ by $h\cdot e^{\alpha}=\<h,\alpha\>e^{\alpha}$ and define $z^h\cdot e^{\alpha}=z^{\<h,\alpha\>}e^{\alpha}$.

Next, we identify $e_{\alpha} $ with $e^{\alpha}$ for $\alpha\in L$ and we set
$$V_L=M(1)\otimes \C[L].$$ The vertex operators $Y(h(-1){\bf 1},z)$ and $Y(e^{\alpha},z)$ associated to $h(-1){\bf 1}$ and $e^{\alpha}$, respectively, are defined as follow:
\begin{eqnarray*}
& &Y(h(-1){\bf 1},z)=h(z)=\sum_{n\in\Z}h(n)z^{-n-1}\\
& &Y(e^{\alpha},z)=\exp\left(\sum_{n=1}^{\infty}\frac{\alpha(-n)}{n}z^n\right)\exp\left(-\sum_{n=1}^{\infty}\frac{\alpha(n)}{n}z^{-n}\right)e_{\alpha}z^{\alpha}.
\end{eqnarray*}
The vertex operator associated with a vector $v=\beta_1(-n_1)...\beta_r(-n_r)e^{\alpha}$ for $\beta_i\in \h$, $n_i\geq 1$ and $\alpha\in L$ is defined by
$$Y(v,z)=:\partial^{(n_1-1)}\beta_1(z)...\partial^{(n_r-1)}\beta_r(z)Y(e^{\alpha},z):,$$
where $\partial^{(n)}=\frac{1}{n!}\left(\frac{d}{dz}\right)^n$ and the normal ordering $:\cdot :$ is an operation which reorders the operators so that $\beta(n)$ ($\beta\in\h,n<0$) and $e_{\alpha}$ to be placed to the left of $\gamma(n)$ ($\gamma\in\h, n\geq 0$) and $z^{\alpha}$.

\begin{thm}\cite{b,flm} $V_L$ is a simple $\Z$-graded vertex algebra with a Virasoro element $\omega=\frac{1}{2}\sum_{a=1}^dh_a(-1)^2\otimes e^0$. Here $\{h_a|1\leq a\leq d\}$ is an orthonormal basis of $\h$. Moreover, $M(1)$ is a $\Z$-graded vertex sub-algebra of $V_L$. \end{thm}

We define a linear automorphism $\theta:V_{L}\rightarrow V_{L}$ by
$$\theta(\b_1(-n_1)\b_2(-n_2)....\b_k(-n_k)e^{\alpha})=(-1)^k\b_1(-n_1)....\b_k(-n_k)e^{-\alpha}.$$ Consequently, $\theta Y(v,z)u=Y(\theta v,z)\theta(u)$ for $u,v\in V_L$. In particular, $\theta$ is an automorphism of $V_L$ and $M(1)$.

For any stable $\theta$-subspace $U$ of $V_{L}$, we denote the $\pm 1$ eigenspace of $U$ for $\theta$ by $U^{\pm}$.
\begin{prop}\cite{dm} $M(1)^+$ and $V_L^+$ are simple vertex algebras.
\end{prop}

Next, we recall constructions of $\theta$-twisted modules of $M(1)$ and $V_L$. We set $\hat{\h}[-1]=\h\otimes t^{\frac{1}{2}}\C[t,t^{-1}]\oplus\C c$. Then $\hat{\h}[-1]$ is an affine Lie algebra with the following commutator formula:
$[\b_1\otimes t^m,\b_2\otimes t^n]=m(\b_1,\b_2)\delta_{m,-n}c\ \ \text{and}\ \ [c,\hat{h}[-1]]=0$ for $\b_i\in\h$, $m,n\in\frac{1}{2}+\Z$. Let $$M(1)(\theta)=U(\hat{\h}[-1])\otimes_{U(\h\otimes t^{1/2}\C[t]\oplus \C c)}\C$$ be the unique irreducible $\hat{\h}[-1]$-module such that $\h\otimes t^{1/2}\C[t]$ acts trivially on $\C$ and $c$ acts on $\C$ as a multiplication by 1.
\begin{prop}\cite{flm} $M(1)(\theta)$ is an irreducible $\theta$-twisted $M(1)$-module.\end{prop}
We will use $\theta$ to denote the automorphism of $\hat{L}$ defined by $\theta(e_{\alpha})=e_{-\alpha}$ and $\theta(\kappa)=\kappa$. We set $K=\{a^{-1}\theta(a)|a\in \hat{L}\}$. For any $\hat{L}/K$-module $T$  such that $\kappa$ acts by the scalar $-1$, we define $$V_L^T=M(1)(\theta)\otimes T.$$ It was shown in \cite{flm} that $V_L^T$ is a $\theta$-twisted $V_L$-module. We define an action of $\theta$ on $M(1)(\theta)$ and $V_L^T$ in the following way:
\begin{eqnarray*}
& &\theta(\b_1(-n_1)\b_2(-n_2)....\b_k(-n_k){\bf 1})=(-1)^k\b_1(-n_1)....\b_k(-n_k){\bf 1},\text{ and }\\
& &\theta(\b_1(-n_1)\b_2(-n_2)....\b_k(-n_k)t)=(-1)^k\b_1(-n_1)....\b_k(-n_k)t).
\end{eqnarray*}
for $\beta_i\in\h$, $n_i\in\frac{1}{2}+\Z_{\geq 0}$ and $t\in T$. We denote by $M(1)(\theta)^{\pm}$ and $V_L^{T,{\pm}}$ the $\pm 1$-eigenspace for $\theta$ of $M(1)(\theta)$ and $V_L^T$, respectively.
\begin{prop}\cite{dlin,dn}\ \ 
\begin{enumerate}
\item $M(1)(\theta)^{\pm}$ are irreducible $M(1)^+$-modules. Furthermore,  
$$\{M(1)^{\pm}, M(1)(\theta)^{\pm}, M(1,\lambda)(\cong M(1,-\lambda))|\lambda\in\h-\{0\}\}$$ is the set of all inequivalent irreducible $M(1)^+$-modules.
\item Let $\chi$ be a central character of $\hat{L}/K$ such that $\chi(\iota(\kappa))=-1$ and $T_{\chi}$ the irreducible $\hat{L}/K$-module with central character $\chi$. Then $V_L^{T_{\chi},\pm}$ are irreducible $V_L^+$-modules.
\end{enumerate}
\end{prop}

\subsection{Zhu algebras $A(M(1)^+)$, $A(V_{\Z\alpha}^+)$}

In this subsection, we recall the Zhu algebras of $M(1)^+$, and $V_{\Z\alpha}^+$ when $\Z\alpha$ is an even lattice such that $\<\alpha,\alpha\>\neq 0$. 

Following the subsection 2.2, we let $L$ be a rank $d$ even lattice with a nondegenerate symmetric $\Z$-bilinear form $\<\cdot,\cdot\>$. We set $\h=\C\otimes_{\Z}L$ and extend $\<\cdot,\cdot\>$ to a $\C$-bilinear form on $\h$. 
Let $\{ h_a|1\leq a\leq d\}$ be an orthonormal basis of $\h$, and set $\omega_a=\omega_{h_a}=\frac{1}{2}h_a(-1)^2{\bf 1}$ and $J_a=h_a(-1)^4{\bf 1}-2h_a(-3)h_a(-1){\bf 1}+\frac{3}{2}h_a(-2)^2{\bf 1}$.  Note that vectors $\omega_a$, and $J_a$ generate a vertex operator algebra $M(1)^+$ associated to the one-dimensional vector space $\C h_a$. 

Following \cite{ad}, we set $S_{ab}(m,n)=h_a(-m)h_b(-n){\bf 1}$, and define $E^u_{ab}$, $E^t_{ab}$, and $\Lambda_{ab}$ as follows:
\begin{eqnarray*}
E^u_{ab}&=&5S_{ab}(1,2)+25S_{ab}(1,3)+36S_{ab}(1,4)+16S_{ab}(1,5)\ \ (a\neq b),\\
E^u_{aa}&=&E^u_{ab}*E^u_{ba},\\
E^t_{ab}&=&-16(3S_{ab}(1,2)+14S_{ab}(1,3)+19S_{ab}(1,4)+8S_{ab}(1,5)\ \ (a\neq b),\\
E^t_{aa}&=&E^t_{ab}*E^t_{ba},\\
\Lambda_{ab}&=&45S_{ab}(1,2)+190S_{ab}(1,3)+240S_{ab}(1,4)+96S_{ab}(1,5).
\end{eqnarray*}
Next, we give a list of relations of some elements in $A(M(1)^+)$ that we will need to use in later sections.
\begin{prop}\label{relam1}\cite{dn} \ \ 
\begin{enumerate} 
\item For any $a,b,c,d$, we have
\begin{enumerate}
\item $[\Lambda_{ab}]=[\Lambda_{ba}]$,
\item $[E^u _{ab}]*[E^u _{cd}]=\delta_{bc}[E^u _{ad}]$, $[E^t _{ab}]*[E^t _{cd}]=\delta_{bc}[E^t_{ad}],$
\item ${[E^u _{ab}]*[E^t _{cd}]}=[E^t _{cd}]*[E^u _{ab}]=0$,
\item ${[\omega_a]*[E^u_{bc}]}=\delta_{ab}[E^u_{bc}]$,
\item ${[E^u_{bc}]*[\omega_a]}=\delta_{ac}[E^u_{bc}]$,
\item ${[\omega_a]*[E^t_{bc}]}=\left(\frac{1}{16}+\frac{1}{2}\delta_{ab}\right)[E^t_{bc}]$,
\item ${[E^t_{bc}]*[\omega_a]}=\left(\frac{1}{16}+\frac{1}{2}\delta_{ac}\right)[E^t_{bc}]$,
\item ${[\omega_a]*[\Lambda_{bc}]}=[\Lambda_{bc}]*[\omega_a]=0$.
\end{enumerate}
\item For any $a,b,c,d$ such that $a\neq b$, we have
\begin{enumerate}
\item ${[\Lambda_{ab}]*[E^u _{cd}]}=[\Lambda_{ab}]*[E^t_{cd}]=[E^u _{cd}]*[\Lambda_{ab}]=[E^t _{cd}]*[\Lambda_{ab}]=0$,
\item ${[\Lambda_{ab}]*[\Lambda_{ab}]}=4[\omega_{a}]*[\omega_b]-\frac{1}{9}([H_a]+[H_b])-([E^u_{aa}]+[E^u_{bb}])-\frac{1}{4}([E^t_{aa}]+[E^t_{bb}])$,
\item $-\frac{2}{9}[H_a]+\frac{2}{9}[H_b]=2[E^u_{aa}]-2[E^u_{bb}]+\frac{1}{4}[E^t_{aa}]-\frac{1}{4}[E^t_{bb}]$,
\item \begin{eqnarray*}
& &-\frac{4}{135}(2[\omega_a]+13)*[H_a]+\frac{4}{135}(2[\omega_b]+13)*[H_b]\\
&=&4([E^u_{aa}]-[E^u_{bb}])+\frac{15}{32}([E^t_{aa}]-[E^t_{bb}]),
\end{eqnarray*}
\item $[\omega_b]*[H_a]=-\frac{2}{15}([\omega_a]-1)*[H_a]+\frac{1}{15}([\omega_b]-1)*[H_b]$.
\end{enumerate}
Here, $H_{a}=J_a+\omega_a-4\omega_a*\omega_a$.
\end{enumerate}
\end{prop}

Let $A^u$, and $A^t$ be the linear subspaces of $A(M(1)^+)$ spanned by $E^u_{ab}$ and $E^t_{ab}$, respectively. Here, $1\leq a,b\leq d$.
\begin{prop}\label{am1gen}\cite{dn} \ \ 
\begin{enumerate}
\item The spaces $A^u$ and $A^t$ are two sided ideals of $A(M(1)^+)$. Moreover, ideals $A^u$, $A^t$, the units $I^u=\sum_{i=1}^d[E^u_{ii}]$ and $I^t=\sum_{i=1}^d[E^t_{ii}]$ of $A^u$, and $A^t$ are independent of the choice of an orthonormal basis.
\item There are algebra isomorphisms between $A^u$ and $\End M(1)^-(0)$, and between $A^t$ and $\End M(1)(\theta)^-$, respectively. In particular, under the basis $$\{h_1(-1){\bf 1},...,h_d(-1){\bf 1}\}$$ of $M(1)^-(0)$, each $[E^u_{ab}]$ corresponds to the matrix element $E_{ab}$ whose $(a,b)$-entry is 1 and zero elsewhere. Similarly, under the basis $$\{h_1(-\frac{1}{2}){\bf 1},...,h_d(-\frac{1}{2}){\bf 1}\}$$ of $M(1)(\theta)^-(0)$, each $[E^t_{ab}]$ corresponds to the matrix element $E_{ab}$ whose $(a,b)$-entry is 1 and zero elsewhere.
\item The Zhu algebra $A(M(1)^+)$ is generated by $[\omega_a]$, $[J_a]$ for $1\leq a\leq d$, $[\Lambda_{ab}]$ for $1\leq a\neq b\leq d$ and $[E^u_{ab}]$, $[E^t_{ab}]$ for $1\leq a,b\leq d$. 
\item The quotient algebra $A(M(1)^+)/(A^t+A^u)$ is commutative. Furthermore, it is generated by the images of $[\omega_a]$, $[J_a]$ for $1\leq a\leq d$ and $[\Lambda_{ab}]$ for $1\leq a\neq b\leq d$.
\end{enumerate}\end{prop}

Next, we recall the Zhu algebra of $V_L^+$ when $L$ is a rank one nondegenerate even lattice. 

Let $\Z\alpha$ be a rank one nondegenerate even lattice, and let $\{h\}$ be an orthonormal basis of $\h=\C\otimes_{\Z} \Z\alpha$. We set $\omega_{\alpha}=\frac{1}{2}h(-1)^2{\bf 1}$, $J_{\alpha}=h(-1)^4{\bf 1}-2h(-3)h(-1){\bf 1}+\frac{3}{2}h(-2)^2{\bf 1}$ and $H_{\alpha}=J_{\alpha}+\omega_{\alpha}-4\omega_{\alpha}*\omega_{\alpha}$. Also, we define $$ E^{\alpha}=e^{\alpha}+e^{-\alpha}\ \ \text{ and }\ \ F^{\alpha}=e^{\alpha}-e^{-\alpha}.$$
\begin{prop}\label{pos}\cite{dn1, dlm1} Suppose that $\Z\alpha$ is a rank one positive definite even lattice. We then have the following.
\begin{enumerate}
\item The Zhu algebra $A(V_{\Z\alpha}^+)$ is generated by $[\omega_{\alpha}]$, $[J_{\alpha}]$ and $[E^{\alpha}]$.
\item If $\<\alpha,\alpha\>=2k\neq 2$, then $A(V_{\Z\alpha}^+)$ is a semisimple, commutative algebra of dimension $k+7$. Moreover, $A(V_{\Z\alpha}^+)$ satisfies the following identity:
\begin{equation}\label{hak}
[H_{\alpha}]*[E^{\alpha}]=\frac{18(8k-3)}{(4k-1)(4k-9)}\left([\omega_{\alpha}]-\frac{k}{4}\right)\left([\omega_{\alpha}]-\frac{3(k-1)}{4(8k-3)}\right)[E^{\alpha}].
\end{equation}
\item If $\<\alpha,\alpha\>=2$, $A(V_{\Z\alpha}^+)$ is a semisimple algebra such that 
\begin{equation}\label{ha2}
[H_{\alpha}]*[E^{\alpha}]+[E^{\alpha}]*[H_{\alpha}]=-12[\omega_{\alpha}]*\left([\omega_{\alpha}]-\frac{1}{4}\right)*[E^{\alpha}].
\end{equation}
\end{enumerate}
\end{prop}
\begin{prop}\label{neggen}\cite{j} Let $\Z\alpha$ be a rank one negative definite even lattice such that $\<\alpha,\alpha\>=-2k$. Then $A(V_{\Z\alpha}^+)$ is a semisimple commutative associative algebra, and $A(V_{\Z\alpha}^+)$ is spanned by $[{\bf 1}]$, $[\omega_{\alpha}]$, $[E^{\alpha}]$, $[\omega_{\alpha}]*[E^{\alpha}]$  with the following identities:
\begin{eqnarray}
& &\left([\omega_{\alpha}]-\frac{1}{16}\right)*\left([\omega_{\alpha}]-\frac{9}{16}\right)=0\label{wa1}\\
& & [E^{2\alpha}]=(1-2k)2^{8k+1}+k2^{8k+6}[\omega_{\alpha}]\\
& &[E^{2\alpha}]*\left(\frac{1+18k}{1+16k}2^{-8k-1}-\frac{k2^{4-8k}}{1+16k}[\omega_{\alpha}]\right)=1
\end{eqnarray} 
Moreover, there is a polynomial $h([\omega_{\alpha}])\in\C[\omega_{\alpha}]$ such that $[J_{\alpha}]=h([\omega_{\alpha}])$.
\end{prop}
\begin{rem} If $\<\alpha,\alpha\>=-2k<0$, then $E^{2\alpha}$ is an invertible element in $A(V_{\Z\alpha}^+)$.\end{rem}
\begin{coro}\label{ja} Let $\Z\alpha$ be a rank one lattice such that $\<\alpha,\alpha\>=-2k$. Then $J_{\alpha}\equiv\left(\frac{9}{128}-\frac{96}{128}\omega_{\alpha}\right)\mod O(V_{\Z\alpha}^+)$.
\end{coro}
\begin{proof} By Proposition \ref{neggen}, we may write $[J_{\alpha}]$ as $a+b[\omega_{\alpha}]$. Next, we let $[J_{\alpha}]$ and $a+b[\omega_{\alpha}]$ act on top levels of $V^{T_1,+}_{\Z\alpha}$ and $V^{T_1,-}_{\Z\alpha}$ (see Table 4). We obtain that 
\begin{eqnarray}
\frac{3}{128}&=&a+\frac{1}{16}b,\\
-\frac{45}{128}&=&a+\frac{9}{16}b.
\end{eqnarray}
By solving the above equations, we then have that $a=\frac{9}{128}$ and $b=-\frac{96}{128}$.
\end{proof}

The following tables provide complete lists of irreducible $M(1)^+$ (respectively, $V_{\Z\alpha}^+$)-modules, and their corresponding irreducible $A(M(1)^+)$(respectively,  $A(V_{\Z{\alpha}}^+)$)-modules which are the top levels of irreducible $M(1)^+$ (respectively, $V_{\Z\alpha}^+$)-modules. In addition, actions of generators of $A(M(1)^+)$ and $A(V_{\Z\alpha}^+)$ on these irreducible modules are given. 

\begin{table}[htdp]\label{table1}
\caption{Actions of $[\omega_a]$, $[J_a]$, $[H_a]$, $[E^u_{ab}]$, $[E^t_{ab}]$  and $[\Lambda_{ab}]$ on top levels of irreducible $M(1)^+$-modules}
\begin{center}
\begin{tabular}{|c|c|c|c|c|c|}
\hline
& $M(1)^+$&$M(1)^-$&$M(1,\lambda)$ ($\lambda\in\h-\{0\}$)&$M(1)(\theta)^+$&$M(1)(\theta)^-$\\
\hline
&${\bf 1}$&$h_c(-1){\bf 1}$&$e^{\lambda}{\bf 1}$&${\bf 1}$&$h_c(-\frac{1}{2}){\bf 1}$\\
\hline
$[\omega_a]$&0&$\delta_{ac}$&$\frac{\<h_a,\lambda\>^2}{2}$& $\frac{1}{16}$&$\frac{1}{16}+\frac{1}{2}\delta_{ac}$\\
\hline
$[J_{a}]$&0&$-6\delta_{ac}$&$\<h_a,\lambda\>^4-\frac{\<h_a,\lambda\>^2}{2}$&$\frac{3}{128}$&$\frac{3}{128}-\frac{3}{8}\delta_{ac}$\\
\hline
$[H_{a}]$&0&$-9\delta_{ac}$&0&$\frac{9}{128}$&$\frac{9}{128}-\frac{9}{8}\delta_{ac}$\\
\hline
$[E^u_{ab}]$&0&$\delta_{bc}h_a(-1){\bf 1}$&0&0&0\\
\hline
$[E^t_{ab}]$&0&0&0&0&$\delta_{bc}h_a(-\frac{1}{2}){\bf1}$\\
\hline
$[\Lambda_{ab}]$&0&0&$\<h_a,\lambda\>\<h_b,\lambda\>$&0&0\\
\hline
\end{tabular}
\end{center}
\label{default}
\end{table}%
\begin{table}[htdp]\label{table2}
\caption{Actions of $[\omega_{\alpha}]$, $[J_{\alpha}]$, $[H_{\alpha}]$ and $[E^{\alpha}]$ on top levels of irreducible $V_{\Z\alpha}^+$-modules when $\<\alpha,\alpha\>=2k$ and $k>1$.}
\begin{center}
\begin{tabular}{|c|c|c|c|c|c|c|}
\hline
&$V_{\Z\alpha}^+$&$V_{\Z\alpha}^-$&$V_{\Z\alpha+\frac{r}{2k}\alpha}$&$V_{\Z\alpha+\frac{\alpha}{2}}^+$&$V_{\Z\alpha+\frac{\alpha}{2}}^-$&$V_{\Z\alpha}^{T_1,+}$\\
& & &\small{($1\leq r\leq k-1$)}& & &  \\
\hline
&${\bf 1}$&$\alpha(-1){\bf 1}$&$e^{\frac{r}{2k}\alpha}$&$e^{\frac{\alpha}{2}}+e^{-\frac{\alpha}{2}}$&$e^{\frac{\alpha}{2}}-e^{-\frac{\alpha}{2}}$& $t_1$\\
\hline
$[\omega_{\alpha}]$&0&1&$r^2/4k$&$k/4$&$k/4$&$1/16$\\
\hline
$[J_{\alpha}]$&0&$-6$&$(\frac{r^2}{2k})^2-\frac{r^2}{4k}$&$\frac{k^4}{4}-\frac{k^2}{4}$&$\frac{k^4}{4}-\frac{k^2}{4}$&$3/128$\\
\hline
$[H_{\alpha}]$&0&$-9$&0&0&0&$9/128$\\
\hline
$[E^{\alpha}]$&0&0&0&1&$-1$&$2^{-2k+1}$\\
\hline
\end{tabular}
\end{center}
\label{default}
\end{table}%
\begin{table}[htdp]
\begin{center}
\begin{tabular}{|c|c|c|c|}
\hline
&$V_{\Z\alpha}^{T_1,-}$&$V_{\Z\alpha}^{T_2,+}$&$V_{\Z\alpha}^{T_2,-}$\\
\hline
&$\alpha (-\frac{1}{2})t_1$& $t_2$& $\alpha (-\frac{1}{2})t_2$\\
\hline
$[\omega_{\alpha}]$&$9/16$&$1/16$&$9/16$\\
\hline
$[J_{\alpha}]$&$-45/128$&$3/128$&$-45/128$\\
\hline
$[H_{\alpha}]$&$-135/128$&$9/128$&$-135/128$\\
\hline
$[E^{\alpha}]$&$-2^{-2k+1}(4k-1)$&$-2^{-2k+1}$&$2^{-2k+1}(4k-1)$\\
\hline
\end{tabular}
\end{center}
\end{table}%
\begin{table}[htdp]\label{table3}
\caption{Actions of $[\omega_{\alpha}]$, $[J_{\alpha}]$, $[H_{\alpha}]$, and $[E^{\alpha}]$ on top levels of irreducible $V_{\Z\alpha}^+$-modules when $\<\alpha,\alpha\>=2$. Note that $c_{\alpha}$ is a square root of $\epsilon(\alpha,\alpha)$.}
\begin{center}
\begin{tabular}{|c|c|c|c|c|c|c|}
\hline
&$V_{\Z\alpha}^+$&$V_{\Z\alpha}^-$& &$V_{\Z\alpha+\frac{\alpha}{2}}^+$&$V_{\Z\alpha+\frac{\alpha}{2}}^-$&$V_{\Z\alpha}^{T_1,+}$\\
\hline
&${\bf 1}$&$\alpha(-1){\bf 1}$&$F^{\alpha}$&$e^{\frac{\alpha}{2}}+c_{\alpha}e^{-\frac{\alpha}{2}}$&$e^{\frac{\alpha}{2}}-c_{\alpha}e^{-\frac{\alpha}{2}}$& $t_1$\\
\hline
$[\omega_{\alpha}]$&0&1& $1$&$1/4$&$1/4$&$1/16$\\
\hline
$[J_{\alpha}]$&0&$-6$&$3$&0&0&$3/128$\\
\hline
$[H_{\alpha}]$&0&$-9$& 0&0&0&$9/128$\\
\hline
$[E^{\alpha}]$&0&$-2F^{\alpha}$& $2\alpha(-1){\bf 1}$&$c^3_{\alpha}$&$-c^3_{\alpha}$&$1/2$\\
\hline
\end{tabular}
\end{center}
\label{default}
\end{table}%
\begin{table}[htdp]
\begin{center}
\begin{tabular}{|c|c|c|c|}
\hline
&$V_{\Z\alpha}^{T_1,-}$&$V_{\Z\alpha}^{T_2,+}$&$V_{\Z\alpha}^{T_2,-}$\\
\hline
&$\alpha (-\frac{1}{2})t_1$& $t_2$& $\alpha (-\frac{1}{2})t_2$\\
\hline
$[\omega_{\alpha}]$&$9/16$&$1/16$&$9/16$\\
\hline
$[J_{\alpha}]$&$-45/128$&$3/128$&$-45/128$\\
\hline
$[H_{\alpha}]$&$-135/128$&$9/128$&$-135/128$\\
\hline
$[E^{\alpha}]$&$-3/2$&$-1/2$&$3/2$\\
\hline
\end{tabular}
\end{center}
\end{table}%
\begin{table}[htdp]\label{table4}
\caption{Actions of $[\omega_{\alpha}]$, $[J_{\alpha}]$, $[E^{\alpha}]$ and $[E^{2\alpha}]$ on top levels of irreducible $V_{\Z\alpha}^+$-modules when $\<\alpha,\alpha\>=-2k$ and $k\geq 1$.}
\begin{center}
\begin{tabular}{|c|c|c|c|c|}
\hline
& $V_{\Z\alpha}^{T_1,+}$&$V_{\Z\alpha}^{T_1,-}$&$V_{\Z\alpha}^{T_2,+}$&$V_{\Z\alpha}^{T_2,-}$\\
\hline
& $t_1$&$\alpha (-\frac{1}{2})t_1$& $t_2$& $\alpha (-\frac{1}{2})t_2$\\
\hline
$[\omega_{\alpha}]$&$1/16$&$9/16$&$1/16$&$9/16$\\
\hline
$[J_{\alpha}]$&$3/128$&$-45/128$&$3/128$&$-45/128$\\
\hline
$[E^{\alpha}]$&$2^{2k+1}$&$(4k+1)2^{2k+1}$&$-2^{2k+1}$&$-(4k+1)2^{2k+1}$\\
\hline
$[E^{2\alpha}]$&$2^{8k+1}$&$(16k+1)2^{8k+1}$&$2^{2k+1}$&$(16k+1)2^{8k+1}$\\
\hline
\end{tabular}
\end{center}
\label{default}
\end{table}%

\section{Classification of irreducible admissible $V_L^+$-modules when $L$ is a negative definite even lattice of a finite rank }
We will study the Zhu algebra of the vertex algebra $V_L^+$ when $L$ is a negative definite even lattice of a finite rank. We will use this information to classify all irreducible admissible $V_L^+$-modules.

In this section, we assume that {\em $L$ is a rank $d$ even lattice with a negative definite symmetric $\Z$-bilinear form $\<\cdot,\cdot\>$}.  Following \cite{ad}, we set $$V_L^+[\alpha]=M(1)^+\otimes E^{\alpha}\oplus M(1)^-\otimes F^{\alpha}$$ and $$A(V_L^+)(\alpha)=(V_L^+[\alpha]+O(V_L^+))/O(V_L^+)$$ for any $\alpha\in L$. Notice that $A(V_L^+)$ is the sum of $A(V_L^+)(\alpha)$ for all $\alpha\in L$. Let $U$ be a $\Z$-graded vertex sub-algebra of $V_L^+$. The identity map induces an algebra homomorphism from $A(U)$ to $A(V_L^+)$. Hence, for $u\in U$, we use $[u]$ to denote $u+O(U)$ and $u+O(V_L^+)$.

We will study spanning sets of $A(V_L^+)(\alpha)$ for any $\alpha\in L$. 
\begin{lem}\label{spav} Let $\alpha\in L$. Then $A(V_L^+)(\alpha)$ is spanned by the vectors $[u]*[E^{\alpha}]*[v]$ for $u,v\in M(1)^+$.\end{lem}
\begin{proof} The proof is as same as the proof in Lemma 5.2 of \cite{ad}.
\end{proof}

\begin{lem}\label{eliminate} For any indices $a,b$, we have $[E^u_{ab}]=0$ and $[\Lambda_{ab}]=0$.
\end{lem}
\begin{proof} Let $\beta\in L$ such that $\<\beta,\beta\>=-2k$. Let $\{h_i|1\leq i\leq d\}$ be an orthonormal basis of $\h$ so that $h_1\in\C\beta$. We set $g([\omega_{1}])=\frac{1+18k}{1+16k}2^{-8k-1}-\frac{k2^{4-8k}}{1+16k}[\omega_{1}]$. Note that $\omega_1=\omega_{\beta}$.

First, we will show that if $a,b\in\{1,...,d\}$ such that $a\neq 1$, then $[E^u_{ab}]=0$. Recall that $[\omega_1]*[E^u_{ab}]=\delta_{1a}[E^u_{ab}]=0$, and $[E^{2\beta}]=(1-2k)2^{8k+1}+k2^{8k+6}[\omega_1]$ (see Proposition \ref{relam1} 1.(d), and Proposition \ref{neggen}). These imply that
$$[E^{2\beta}]*[E^u_{ab}]=(1-2k)2^{8k+1}[E^u_{ab}].$$ Since $g([\omega_{1}])*[E^{2\beta}]=1$, we then have
\begin{eqnarray*}
[E^{u}_{ab}]&=&g([\omega_{1}])*[E^{2\beta}]*[E^u_{ab}]\\
&=&\frac{1+18k}{1+16k}2^{-8k-1}(1-2k)2^{8k+1}[E^u_{ab}].\\
\end{eqnarray*}
Consequently, $[E^u_{ab}]=0$. Next, we will show that $[E^u_{1b}]=0$ for $b\in\{1,...,d\}$. Since $[\omega_{1}]*[E^u_{1b}]=[E^u_{1b}]$ and 
\begin{eqnarray*}
[E^u_{1b}]&=&\left(g([\omega_{\beta}])*[E^{2\beta}]\right)*[E^u_{1b}]\\
&=&\left(1-\frac{420k^2}{1+16k}\right)[E^u_{1b}],
\end{eqnarray*}
it follows immediately that $[E^u_{1b}]=0$. Hence, for any indices $a,b$, $[E^u_{ab}]=0$.

To show that $[\Lambda_{ab}]=0$ for any indices $a,b$, we only need to use the fact that $[\omega_c]*[\Lambda_{ab}]=0$ for any $c\in\{1,....,d\}$ (see Proposition \ref{relam1} 1.(h)) and follow the first part of the proof of $[E^u_{ab}]=0$ step by step. 

\end{proof}
\begin{lem}\label{etea} For $\alpha\in L-\{0\}$, we let $\{h_1,...,h_d\}$ be an orthonormal basis of $\h$  such that $h_1\in \C\alpha$. Then 
\begin{enumerate}
\item $[E^t_{ab}]*[E^{\alpha}]=[E^{\alpha}]*[E^t_{ab}]$ if $a\neq 1$ and $b\neq 1$.
\item $[E^t_{1b}]*[E^{\alpha}]=-\frac{1}{2\<\alpha,\alpha\>-1}[E^{\alpha}]*[E^t_{1b}]$ if $b\neq 1$.
\item $[E^t_{b1}]*[E^{\alpha}]=-(2\<\alpha,\alpha\>-1)[E^{\alpha}]*[E^t_{b1}]$ if $b\neq 1$.
\item $[E^t_{aa}]*[E^{\alpha}]=[E^{\alpha}]*[E^t_{aa}]$ for $a\in\{1,...,d\}$.
\end{enumerate}
\end{lem}
\begin{proof} 1. is clear. 2., and 3. can be easily obtained by using Lemma \ref{eliminate} combines with the proof of Lemma 7.5 of \cite{ad}. 

Next, we will prove 4. Recall that $E_{aa}^t=E_{ab}^t*E^t_{ba}$ where $b\neq a$. If $a\neq 1$, then, by 1., we have that $[E^t_{aa}]*[E^{\alpha}]=(E_{ab}^t*E^t_{ba})*[E^{\alpha}]=[E^{\alpha}]*[E^t_{aa}].$ Similarly, by 2., and 3., we then have that $[E^t_{11}]*[E^{\alpha}]=([E^t_{1b}]*[E^t_{b1}])*[E^{\alpha}]=[E^{\alpha}]*([E^t_{1b}]*[E^t_{b1}])=[E^{\alpha}]*[E^t_{11}]$.\end{proof}

\begin{prop} \label{it} Let $I^t$ be the unit of the simple algebra $A^t$. Then for any $\alpha\in L$, $I^t*[E^{\alpha}]=[E^{\alpha}]*I^t$.
\end{prop}
\begin{proof} It follows immediately from Lemma \ref{etea}.
\end{proof}
\begin{prop}\label{spanavla} For any $\alpha\in L$, we have 
\begin{eqnarray*}
A(V_L^+)(\alpha)&=&span_{\C}\{[u]*[E^{\alpha}]| u\in M(1)^+\}\\
&=&span_{\C}\{[E^{\alpha}]*[u]| u\in M(1)^+\}.\end{eqnarray*}\end{prop}
\begin{proof} Let $\alpha\in L-\{0\}$ and let $\{h_1,...,h_d\}$ be an orthonormal basis of $\h$ such that $h_1\in\C\alpha$. We will show that for $v\in M(1)^+$, $$[E^{\alpha}]*[v]\in span_{\C}\{[u]*[E^{\alpha}]| u\in M(1)^+\}.$$ Since $\{[v]+O(V_L^+)|v\in M(1)^+\}$ is generated by $[\omega_a]$, $[H_a]$, $[E^t_{ab}]$ where $1\leq a,b\leq d$ (see Proposition \ref{am1gen} and Lemma \ref{eliminate}), it is enough to show that for $a,b\in\{1,...,d\}$, 
$$[E^{\alpha}]*[\omega_a],[E^{\alpha}]*[H_a],[E^{\alpha}]*[E^t_{ab}]\in span_{\C}\{[u]*[E^{\alpha}]| u\in M(1)^+\}.$$ 
By Lemma \ref{etea}, we have that for $a,b\in\{1,...,d\}$,$$[E^{\alpha}]*[E^{t}_{ab}]\in span_{\C}\{[u]*[E^{\alpha}]| u\in M(1)^+\}.$$ Clearly, for $a\in \{1,...,d\}$, $E^{\alpha}*[\omega_a]=[\omega_a]*[E^{\alpha}]$. Since $[J_{1}]\in\C[\omega_1]$ and $[H_1]=[J_1]+[\omega_1]-4[\omega_1]*[\omega_1]$, it follows that $[E^{\alpha}]*[H_1]=[H_1]*[E^{\alpha}]$. Furthermore, by  Proposition \ref{relam1} 2.(c) and Lemma \ref{eliminate}, we have that for $a\in\{2,...,d\}$,
$$[H_a]=[H_1]-\frac{9}{8}[E^t_{aa}]+\frac{9}{8}[E^t_{11}].$$
This implies that $[E^{\alpha}]*[H_a]\in span_{\C}\{[u]*[E^{\alpha}]| u\in M(1)^+\}.$

The proof of $A(V_L^+)(\alpha)=span_{\C}\{[E^{\alpha}]*[u]| u\in M(1)^+\}$ is very similar to the above.
\end{proof}
\begin{coro}\label{hao} Let $\alpha\in L-\{0\}$. We set $\{h_1,...,h_d\}$ be an orthonormal basis of $\h$ such that $h_1\in\C\alpha$. Then $$[H_a]=[H_1]-\frac{9}{8}[E^t_{aa}]+\frac{9}{8}[E^t_{11}]\text{\ \ \ \ \ \   for }a\in\{2,...,d\} .$$
\end{coro}

Next, we will classify all irreducible $A(V_L^+)$-modules.
\begin{de} Let $W$ be an $A(V_L^+)$-module. An element $u\in A(V_L^+)$ is called semisimple on $W$ if $u$ acts diagonally on $W$.
\end{de}

\begin{lem}\label{wsemi} Let $W$ be an irreducible $A(V_L^+)$-module such that $A^tW=0$. Then every element in $A(M(1)^+)$ is semisimple on $W$.\end{lem}
\begin{proof} Recall that for any $\alpha\in L-\{0\}$, $A(V_{\Z\alpha}^+)$ is a semisimple commutative algebra and $[\omega_{\alpha}]$, $[J_{\alpha}]$ act semisimply on $W$. By Lemma \ref{eliminate}, we can conclude further that $A^u$ and $[\Lambda_{ab}]$ acts as zero on $M$ for all $1\leq a,b\leq d$.  Since $A(M(1)^+)$ is generated by $[\omega_a]$, $[J_a]$ for $1\leq a\leq d$, $[\Lambda_{ab}]$ for $1\leq a\neq b\leq d$ and $A^u$, $A^t$, we then have that every element in $A(M(1)^+)$ is semisimple on $W$. 
\end{proof}
\begin{lem}\label{wcon} Let $W$ be an irreducible $A(V_L^+)$-module such that $A^tW=0$. Then any element in $A(M(1)^+)$ acts as a constant on $W$. Consequently, for $\alpha\in L$, the action of $[E^{\alpha}]$ and the action of every element in $A(M(1)^+)$ commute on $W$.
\end{lem}
\begin{proof} 
By Proposition \ref{am1gen} 4. and Lemma \ref{wsemi}, we can conclude that $W$ is a direct sum of one-dimensional irreducible $A(M(1)^+)$-modules on which $[\Lambda_{ab}]$, $A^u$ and $A^t$  act as zero on $W$. Here $1\leq a, b\leq d$. This implies that an irreducible $A(M(1)^+)$-submodule of $W$ is isomorphic to either $M(1)^+(0)$, or $M(1)(\theta)^+(0)$. 

Since $([\omega_{\alpha}]-\frac{1}{16})*([\omega_{\alpha}]-\frac{9}{16})\equiv 0\mod O(V_{\Z\alpha}^+)$ for every $\alpha\in L-\{0\}$ (see Proposition \ref{neggen}), it follows that $[\omega_{\alpha}]$ acts on $W$ either as $1/16$ or $9/16$. By Corollary \ref{ja} and the fact that $[H_{\alpha}]=[J_{\alpha}]+[\omega_{\alpha}]-4[\omega_{\alpha}]*[\omega_{\alpha}]$, we have that $[J_{\alpha}]$ acts either as $\frac{3}{138}$ or $-\frac{45}{128}$ on $W$ and $[H_{\alpha}]$ acts either as $\frac{9}{128}$ or $-\frac{35}{128}$ on $W$. Let $\{h_1,...,h_d\}$ be an orthonormal basis so that $h_1\in\C\alpha$. Then by Corollary \ref{hao} and the fact that $A^tW=0$, we can conclude that, on $W$, 
$$[H_a]=[H_1]=[H_{\alpha}]$$ for all $a\in\{1,...,d\}$. By Proposition \ref{relam1} 2.(d), 2.(e), we have that, on $W$, $[\omega_{a}]=\frac{1}{16}$ for all $a\in\{1,...,d\}$. Consequently, for $a\in\{1,...,d\}$
$$[H_a]\text{ acts as }\frac{9}{128}.$$
By Table 1, we have that an irreducible $A(M(1)^+)$ submodule of $W$ is isomorphic to $M(1)(\theta)^+(0)$. Since $[\omega_a]$, $[H_a]$ act on $W$ as $1/16$, and $9/128$, respectively for all $a\in\{1,...,d\}$, this implies that every element in $A(M(1)^+)$ acts as a constant on $W$ and for any $\alpha\in L$, the action of $[E^{\alpha}]$ and the action of every element in $A(M(1)^+)$ commute on $W$.
\end{proof}
\begin{coro} $W$ is the direct sum of $M(1)(\theta)^+(0)$.
\end{coro}

Let $W$ be an irreducible $A(V_L^+)$-module such that $A^tW=0$. We will show that $W$ is isomorphic to an irreducible $\hat{L}/K$-module $T_{\chi}$ for some $\chi$.

Let $\alpha\in L$ such that $\<\alpha,\alpha\>=-2k$.
Recall that $$E^{2\alpha}=(1-2k)2^{8k+1}+k2^{8k+6}[\omega_{\alpha}]$$ and $[\omega_{\alpha}]$ acts as $\frac{1}{16}$ on $W$. Hence $E^{2\alpha}$ acts as $2^{-4\<\alpha,\alpha\>+1}$ on $W$ and on $V_L^{T_{\chi}}(0)$ for any $\chi$.

For $\alpha\in L-\{0\}$, we set $$B_{\alpha}=2^{\<\alpha,\alpha\>-1}E^{\alpha}.$$ Also, we let $B_0={\bf 1}$. Since $E^{2\alpha}$ acts as $2^{-4\<\alpha,\alpha\>+1}$ on $W$ and on $V_L^{T_{\chi}}(0)$ for any $\chi$, we can conclude that $B_{2\alpha}=1$ on $W$, and on $V_L^{T_{\chi}}(0)$ for any $\chi$. 
\begin{lem}\label{baction} For $\alpha,\beta\in L$, $[B_{\alpha}]*[B_{\beta}]=\epsilon (\alpha,\beta)[B_{\alpha+\beta}]$ on $W$. \end{lem}

\begin{proof} We will follow the idea in \cite{ad}. Let $\alpha,\beta\in L$.
\begin{description}
\item[case 1:] $\<\alpha,\beta\> < 0$. 

Then $E^{\alpha}*E^{\beta}\in V_L^+[\alpha+\beta]$. By Proposition \ref{spanavla}, we have that 
$$[B_{\alpha}]*[B_{\alpha}]=[u]*[B_{\alpha+\beta}]$$ where $u\in M(1)^+$. Since $[u]$ acts on $M(1)(\theta)^+(0)$ as a constant, say $p$, we obtain 
$$[B_{\alpha}]*[B_{\beta}]=p[B_{\alpha+\beta}]$$ on $W$. On the other hand,
$[B_{\alpha}]*[B_{\beta}]=\epsilon(\alpha,\beta)[B_{\alpha+\beta}]$ on $V_L^{T_{\chi,+}}(0)$ for any $\chi$. Hence, $[B_{\alpha}]*[B_{\beta}]=\epsilon (\alpha, \beta)[B_{\alpha+\beta}]$ on $W$.

\item[case 2:] $\<\alpha,\beta\> > 0$. 

Note that 
\begin{equation}\label{abt}[B_{\alpha}]*[B_{\beta}]=[v]*[B_{-\alpha+\beta}]\end{equation}
for some $v\in M(1)^+$ since $E^{\alpha}*E^{\beta}\in V_L^+[-\alpha+\beta]$.
\begin{description}
\item[case 2.1:] $\<\beta,\beta\><\<\alpha,\alpha\>$. 

Then $\<\alpha+\beta, -2\alpha+2\beta\> <0$, and we have
\begin{eqnarray}
[B_{\alpha+\beta}]&=&[B_{\alpha+\beta}]*[B_{-2\alpha+2\beta}]\nonumber\\
&=&\epsilon(\alpha+\beta,2\alpha-2\beta)[B_{-\alpha+3\beta}]\\
&=&[B_{-\alpha+3\beta}]\label{e2}
\end{eqnarray}
on $W$ and on $V_L^{T_{\chi},+}(0)$ for any $\chi$.
This implies 
\begin{eqnarray}
[B_{\alpha}]*[B_{\beta}]&=&[v]*[B_{-\alpha+\beta}]\nonumber\\
&=&[v]*([B_{2\beta}]*[B_{-\alpha+\beta}])\nonumber\\
&=&[v]*\epsilon(2\beta, -\alpha+\beta)[B_{-\alpha+3\beta}]\nonumber\\
&=&[v]*[B_{\alpha+\beta}]\text{ (by (\ref{e2}))}\nonumber
\end{eqnarray}
on $W$ and on $V_L^{T_{\chi},+}(0)$ for any $\chi$.
Next, by following the proof of {\bf case 1}, we then have
$$[B_{\alpha}]*[B_{\beta}]=\epsilon(\alpha,\beta)[B_{\alpha+\beta}]\text{ on } W.$$

\item[case 2.2] $\<\alpha,\alpha\> < \<\beta,\beta\>$.  

Then $\<\alpha+\beta,-2\beta+2\alpha\><0$. So, we have
\begin{eqnarray}
[B_{\alpha+\beta}]&=&[B_{\alpha+\beta}]*[B_{2\alpha-2\beta}]\nonumber\\
&=&\epsilon(\alpha+\beta,-2\beta+2\alpha)[B_{3\alpha-\beta}]\\
&=&[B_{3\alpha-\beta}]\label{e3}
\end{eqnarray}
on $W$ and on $V_L^{T_{\chi},+}(0)$ for any $\chi$. Consequently,
\begin{eqnarray}
[B_{\alpha}]*[B_{\beta}]&=&[v]*[B_{-\alpha+\beta}]\nonumber\\
&=&[v]*[B_{\alpha-\beta}]\nonumber\\
&=&[v]*([B_{2\alpha}])*[B_{\alpha-\beta}])\nonumber\\
&=&[v]*[B_{3\alpha-\beta}]\nonumber\\
&=&[v]*[B_{\alpha+\beta}] \text{ (by (\ref{e3}))}.
\end{eqnarray}
As in {\bf case 1}, we have that on $W$,
$$[B_{\alpha}]*[B_{\beta}]=\epsilon(\alpha,\beta)[B_{\alpha+\beta}].$$

\item[case 2.3:] $\<\alpha,\alpha\>=\<\beta,\beta\>$. 


Since $B_{-2\beta+2\alpha}*B_{\alpha+\beta}\in V_L^+[3\alpha-\beta]+V_L^+[-3\beta+\alpha]$, we then have 
\begin{eqnarray}
[B_{\alpha+\beta}]&=&[B_{-2\beta+2\alpha}]*[B_{\alpha+\beta}]\nonumber\\
&=&[u]*[B_{3\alpha-\beta}]+[w]*[B_{-3\beta+\alpha}]\text{ where }u,w\in M(1)^+ \nonumber\\
&=&[u]*([B_{2\alpha}]*[B_{\alpha-\beta}])+[w]*([B_{-2\beta}]*[B_{-\beta+\alpha}]) \nonumber\\
&=&([u+v])*[B_{\alpha-\beta}] \nonumber\\
&=&c[B_{\alpha-\beta}]\text{ for some constant }c\nonumber
\end{eqnarray}
on $W$ and on $V_L^{T_{\chi},+}(0)$ for any $\chi$. On the other hand, $[B_{-2\beta}]*[B_{\alpha+\beta}]=\epsilon(-2\beta,\alpha+\beta)[B_{\alpha-\beta}]=[B_{\alpha-\beta}]$ on $V_L^{T_{\chi},+}(0)$ for any $\chi$. Hence, $c=1$ and 
\begin{equation}\label{e4}[B_{\alpha+\beta}]=[B_{\alpha-\beta}].\end{equation}
By (\ref{abt}) and (\ref{e4}), we have that on $W$ and on $V_L^{T_{\chi},+}(0)$ for any $\chi$, $[B_{\alpha}]*[B_{\beta}]=[v]*[B_{\alpha+\beta}]$.
By following the proof of the {\bf case 1}, we then obtain that 
$$[B_{\alpha}]*[B_{\beta}]=\epsilon(\alpha,\beta)[B_{\alpha+\beta}].$$
\end{description}

\item[case 3:] $\<\alpha,\beta\>=0$. 

Then $[E^{\alpha}]*[E^{\beta}]=\epsilon(\alpha,\beta)([E^{\alpha+\beta}]+[E^{\alpha-\beta}])$. To show that $[B_{\alpha}]*[B_{\beta}]=\epsilon(\alpha,\beta)[B_{\alpha+\beta}]$, we only need to modify the proof of the {\bf case 2}. 
\end{description}
\end{proof}

\begin{thm}\label{cam1} Let $W$ be an irreducible $A(V_L^+)$-module such that $A^tW=0$. Then there exists an irreducible $\hat{L}/K$-module $T_{\chi}$ with central character $\chi$ such that $W\cong T_{\chi}=V_L^{T_\chi, +}(0)$. Here, $\chi(\iota(\kappa))=-1$.\end{thm}
\begin{proof} We will follow the proof of the proposition 6.6 in \cite{ad}. First, we construct a map from $\hat{L}$  to $GL(W)$ by sending $e_{\alpha}$ to $B_{\alpha}$  and $\kappa$ to $-1$. This map will induce a group homomorphism from $\hat{L}/K$ to $GL(W)$ since $\theta(e_{\alpha})=e_{-\alpha}$ and $B_{\alpha}=B_{-\alpha}$. Since the action of $A(M(1)^+)$ on $W$ commute with the action of $[B_{\alpha}]$ for all $\alpha\in L$, by Proposition \ref{spanavla}, we can conclude that $W$ is an irreducible $\hat{L}/K$-module. Hence, $W$ is isomorphic to $T_{\chi}=V_L^{T_\chi, +}(0)$ for some central character $\chi$.\end{proof}

Next, we let $M$ be an irreducible $A(V_L^+)$-module such that $A^t M \neq 0$. Since $A^t$ is a simple algebra, it implies that $M$ contains a simple $A^t$-module which is isomorphic to $M(1)(\theta)^-(0)$. 
\begin{lem}\label{am+act} $M=A^tM$ and $M$ is a direct sum of $M(1)(\theta)^-(0)$.\end{lem}
\begin{proof} We will follow the proof in \cite{ad}. We set $M^0=\{m\in M|A^tm=0\}$. Since $A^t$ is a two-sided ideal of $A(M(1)^+)$, it follows immediately that $M^0$ is an $A(M(1)^+)$-module. We will show that $M^0$ is an $A(V_L^+)$-module. Let $m\in M^0$ and $u\in A^t$. By Proposition \ref{it}, we have $$[u]\cdot([E^{\alpha}]m)=([u]*[I^t])\cdot([E^{\alpha}]m)=([u]*[E^{\alpha}])\cdot(I^tm)=0\text{ for all }\alpha\in L.$$ By Proposition \ref{spanavla}, we can conclude that $M^0$ is an $A(V_L^+)$-submodule of $M$. Since $M$ is an irreducible $A(V_L^+)$-module such that $A^tM\neq 0$, it implies that $M^0=0$. Hence, $M=A^tM$, and $M$ is a direct sum of $M(1)(\theta)^-(0)$.
\end{proof}
 
Next, for any $\alpha\in L-\{0\}$, we define 
\begin{equation}
\tilde{B}_{\alpha}=2^{\<\alpha,\alpha\>-1}(E^{\alpha}-\frac{2\<\alpha,\alpha\>}{2\<\alpha,\alpha\>-1}E^t_{11}*E^{\alpha})
\end{equation}
and we set $\tilde{B}_0=1$. Here, $E^t_{11}$ is defined with respect to an orthonormal basis $\{h_a|1\leq a\leq d\}$ of $\h$ such that $h_1\in \C\alpha$. Note that for $\alpha,\beta\in L$, $[\tilde{B}_{\alpha}]*[\tilde{B}_{\beta}]=\epsilon(\alpha,\beta)[\tilde{B}_{\alpha+\beta}]$ on $V_L^{T_{\chi},-}(0)$ for any $\chi$. Also, by using Table 1 and the fact that $[E^{2\alpha}]=(1-2k)2^{8k+1}+k2^{8k+6}[\omega_{\alpha}]$ where $\<\alpha,\alpha\>=-2k$, one can show that for $\alpha\in L$, 
$$[\tilde{B}_{2\alpha}]=1\text{ on }M.$$
We will prove that $M$ is a $\hat{L}/K$-module.
\begin{lem}\label{comm} For any $1\leq a,b\leq d$, $[\tilde{B}_{\alpha}]$ and $[E^t_{ab}]$ commute on $M$. Consequently, $[\tilde{B}_{\alpha}]$ commutes with the action of $[u]$ where $u$ is any member of $M(1)^+$.
\end{lem}
\begin{proof} By using Proposition \ref{relam1} 1.(b), and Lemma \ref{etea}, we can directly show that $[\tilde{B}_{\alpha}]$ commutes with $[E^t_{ab}]$ for $1\leq a,b\leq d$. Moreover, by Lemma \ref{am+act}, we can conclude further that the action of $[\tilde{B}_{\alpha}]$ on $M$ commutes with the action of $[u]$ on $M$. Here $u$ is any member of $M(1)^+$.
\end{proof}
\begin{lem}\label{modl} For $\alpha,\beta\in L$, we have 
\begin{equation}\label{am}[\tilde{B}_{\alpha}]*[\tilde{B}_{\beta}]=\epsilon(\alpha,\beta)[\tilde{B}_{\alpha+\beta}]\text{ on } M.\end{equation}
\end{lem}
\begin{proof} Let $\alpha,\beta\in L$ such that $\<\alpha,\beta\><0$. Using Proposition \ref{spanavla}, Lemma \ref{am+act}, Lemma \ref{comm}, and the facts that $[E^{\alpha}]*[E^{\beta}]\in V_L^+[\alpha+\beta]$, we then have that 
$$[\tilde{B}_{\alpha}]*[\tilde{B}_{\beta}]=[u]*[\tilde{B}_{\alpha+\beta}]\text{ where }u\in A^t.$$ Since $[\tilde{B}_{\alpha}]*[\tilde{B}_{\beta}]=\epsilon(\alpha,\beta)[\tilde{B}_{\alpha+\beta}]$ on any irreducible $A(V_L^+)$-modules $V_L^{T_{\chi},-}(0)$ and $M$ is a direct sum of $M(1)(\theta)^-(0)$, we can conclude that $[u]=\epsilon(\alpha,\beta)I^t$. Hence, $[\tilde{B}_{\alpha}]*[\tilde{B}_{\beta}]=\epsilon(\alpha,\beta)[\tilde{B}_{\alpha+\beta}]$ on $M$. The proof of an equation (\ref{am}) for the case $\alpha,\beta\in L$ such that $\<\alpha,\beta\>\geq 0$ is similar to the proof of Lemma \ref{baction}.
\end{proof}

\begin{thm}\label{thm2} Let $M$ be an irreducible $A(V_L^+)$-module such that $A^tM\neq 0$. Then there exists an irreducible $\hat{L}/K$-module $T_{\chi}$ with central character $\chi$ such that $M\cong \h(-1/2)\otimes T_{\chi}=V_L^{T_{\chi}, -}(0)$. Here, $\chi(\iota(\kappa))=-1$.
\end{thm}
\begin{proof} By following the proof of Theorem \ref{cam1}, we can show that $M$ is a $\hat{L}/K$-module. Moreover, by Lemma \ref{comm}, we can conclude that $M=\h(-\frac{1}{2})\otimes U$  where $U$ is a $\hat{L}/K$-module. Since $M$ is an irreducible $A(V_L^+)$-module, it follows that $U$ is an irreducible $\hat{L}/K$-module.
\end{proof}
The following is a consequent of Theorems \ref{zhu}, \ref{cam1}, and \ref{thm2}.
\begin{prop} Let $L$ be a negative definite even lattice. Then the set of all inequivalent irreducible admissible $V_L^+$-module is 
\begin{eqnarray*}
\{\ \ V_L^{T_{\chi},\pm}&|& T_{\chi}\text{ is irreducible }\hat{L}/K\text{-module with central character }\chi\\
& &\text{ such that }\chi(\iota(\kappa))=-1\}.\end{eqnarray*}
\end{prop}

\section{Classification of irreducible $V_L^+$-modules when $L$ is a nondegenerate even lattice}
We will classify all irreducible admissible $V_L^+$-modules when $L$ is an even lattice of a finite rank equipped with a symmetric nondegenerate $\Q$-valued $\Z$-bilinear form that is neither positive definite nor negative definite. In particular, we will show that all irreducible admissible $V_L^+$-modules are $V_L^{T_{\chi},\pm}$ where $T_{\chi}$ are irreducible $\hat{L}/K$-modules with central character $\chi$. Here, $\chi(\iota(\kappa))=-1$.

For this section, we assume that {\em $L$ is a rank $d$ even lattice equipped with a symmetric nondegenerate $\Q$-valued $\Z$-bilinear form that is neither positive definite nor negative definite.} First, we will recall several propositions that will be useful later.
\begin{prop}\label{tb}\cite{la} Let $k$ be an ordered field and let $E$ be a nonzero finite dimensional vector space over $k$ with a non-degenerate symmetric form $\<\ \ ,\ \ \>$. Then $E$ has an orthogonal basis. Moreover, there exists an integer $r\geq 0$ such that if $\{ v_1,...,v_n\}$ is an orthogonal basis of $E$ then precisely $r$ among the $n$ elements $\<v_1,v_1\>,....,\<v_n,v_n\>$ are $>0$ and $n-r$ among these elements are $<0$. If we assume that every positive element of $k$ is a square. Then there exists an orthogonal basis $\{v_1,...,v_n\}$ of $E$ such that $\<v_i,v_i\>=1$ for $i\leq r$ and $\<v_i,v_i\>=-1$ for $i>r$ and $r$ is uniquely determined.\end{prop}
\begin{prop}\label{spm1}\cite{ad} Let $\alpha\in L-\{0\}$. Then $V_L^+[\alpha]$ is spanned by the vectors $u*E^{\alpha}$ and $u*h(-n)F^{\alpha}$ for $u\in M(1)^+$, $h\in\h$, and $n\geq 1$.
\end{prop}

Next, we will study the spanning set of $A(V_L^+)$. The following proposition resembles results in section 3. In fact the proof of this proposition is the same as the proof of the results in section 3.
\begin{prop}\label{etea2}\ \
\begin{enumerate}
\item For $\alpha\in L$ such that $\<\alpha,\alpha\>\neq 0$, $A(V_L^+)(\alpha)$ is spanned by the vectors $[u]*[E^{\alpha}]*[v]$ for $u,v\in M(1)^+$.
\item For any indices $a,b$ $[E^u_{ab}]=0$ and $[\Lambda_{ab}]=0$.
\item For $\alpha\in L$ such that $\<\alpha,\alpha\>\neq 0$, we set $\{h_1,...,h_d\}$ be an orthonormal basis of $\h$ such that $h_1\in \C\alpha$. Then 
\begin{enumerate}
\item $[E^t_{ab}]*[E^{\alpha}]=[E^{\alpha}]*[E^t_{ab}]$ if $a\neq 1$ and $b\neq 1$.
\item $[E^t_{1b}]*[E^{\alpha}]=-\frac{1}{2\<\alpha,\alpha\>-1}[E^{\alpha}]*[E^t_{1b}]$ if $b\neq 1$.
\item $[E^t_{b1}]*[E^{\alpha}]=-(2\<\alpha,\alpha\>-1)[E^{\alpha}]*[E^t_{b1}]$ if $b\neq 1$.
\item $[E^t_{aa}]*[E^{\alpha}]=[E^{\alpha}]*[E^t_{aa}]$ for $a\in \{1,...,d\}$.
\end{enumerate}
\item Let $I^t$ be the unit of the simple algebra $A^t$. Then for any $\alpha\in L$ such that $\<\alpha,\alpha\>\neq 0$, $[I^t]*[E^{\alpha}]=[E^{\alpha}]*[I^t]$.
\item For $\alpha\in L-\{0\}$, we set $\{h_1,...,h_d\}$ be an orthonormal basis of $\h$ so that $h_1\in \C{\alpha}$. Then $[H_a]=[H_1]-\frac{9}{8}[E^t_{aa}]+\frac{9}{8}[E^t_{11}]$ for $a\in\{2,...,d\}$.
\end{enumerate}
\end{prop}

\begin{thm}\label{spangen} \ \ 
\begin{enumerate}
\item For any $\alpha\in L$ such that $\<\alpha,\alpha\>\neq 0$, we have 
\begin{eqnarray*}
A(V_L^+)(\alpha)&=&span_{\C}\{[u]*[E^{\alpha}]|u\in M(1)^+\}\\
&=&span_{\C}\{[E^{\alpha}]*[u]|u\in M(1)^+\}.
\end{eqnarray*}
\item Let $\alpha\in L$ such that $\<\alpha,\alpha\>=0$. Then there exists $\beta\in L$ such that $\<\beta,\beta\><0$, $\<\alpha,\beta\><0$, and $A(V_L^+)(\alpha)\subset A(V_L^+)(\alpha+2\beta)$. 
\item $A(V_L^+)$ is spanned by $A(V_L^+)(\alpha)$ for all $\alpha\in L$ such that $\<\alpha,\alpha\>\neq 0$.
\end{enumerate}\end{thm}
\begin{proof}  For 1., the proof is similar to the proof of Proposition \ref{spanavla}.


For 2., we let $\alpha\in L-\{0\}$ such that $\<\alpha,\alpha\>=0$. By Proposition \ref{tb}, we can conclude that there exist $\gamma, \lambda\in L_{\Q}$ such that $\alpha=\gamma+\lambda$, $\<\gamma,\gamma\>>0$, $\<\lambda,\lambda\><0$ and $\<\gamma,\lambda\>=0$. Since $\lambda\in{L_{\Q}}$, we know that there exists a positive integer $m$ such that $m\lambda\in L$ and $\<m\lambda,m\lambda\>$ is an even negative integer. For convenience, we set $\beta=m\lambda$. Clearly, $\<\alpha+2\beta,\alpha+2\beta\>$ is an even negative integer. Since
$$[E^{\alpha}]=(g([\omega_{\beta}])*[E^{2\beta}])*[E^{\alpha}],\ \ [h(-n)F^{\alpha}]=(g([\omega_{\beta}])*[E^{2\beta}])*[h(-n)F^{\alpha}],$$ and $$[E^{2\beta}]*[E^{\alpha}],\ \  [E^{2\beta}]*[h(-n)F^{\alpha}]\in V_L^+[\alpha+2\beta]\text{ for all }h\in\h, n\geq 1,$$ by Proposition \ref{spm1} we can conclude that $A(V_L^+)(\alpha)\subset A(V_L^+)(\alpha+2\beta)$. 

3. is clear.
\end{proof}
\begin{rem} Let $\alpha\in L$ such that $\<\alpha,\alpha\>=0$ and let $\rho\in L$ such that $\<\rho,\rho\><0$ and $\<\rho,\alpha\><0$. Then $A(V_L^+)(\alpha)\subset A(V_L^+)(\alpha+2\rho)$.
\end{rem}

Next, we will classify all irreducible admissible $V_L^+$-modules. We let $W$ be an irreducible $A(V_L^+)$-module such that $A^tW=0$.
\begin{lem} Then $W$ is a direct sum of $M(1)(\theta)^+(0)$ and every element in $A(M(1)^+)$ acts as a constant on $W$. Consequently, $[E^{\alpha}]$ commutes with the action of $[u]$ where $u$ is any member of $M(1)^+$.
\end{lem}
\begin{proof} The proof is very similar to Lemma \ref{wsemi}, and Lemma \ref{wcon}.
\end{proof} 
Next, we will show that $W$ is isomorphic to $V_L^{T_{\chi},+}(0)$ for some central character $\chi$ such that $\chi(\iota(\kappa))=-1$. Following section 3, we set 
$$B_{\alpha}=2^{\<\alpha,\alpha\>-1}E^{\alpha}\text{ for } \alpha\in L-\{0\}$$ and $B_{0}=1$. Note that if $\alpha\in L$ such that $\<\alpha,\alpha\><0$, then $[B_{2\alpha}]={\bf 1}$ on $W$ and on $V_L^{T_{\chi}}(0)$ for all central character $\chi$ such that $\chi(\iota(\kappa))=-1$.
\begin{lem}\label{lb2action} Let $\alpha\in L$ such that $\<\alpha,\alpha\>=0$. Then there exists $\beta\in L$ such that $\<\beta,\beta\><0$, $\<\alpha,\beta\><0$ and $[B_{\alpha}]=[B_{\alpha+2\beta}]$ on $W$, and on $V_L^{T_{\chi},+}(0)$ for any $\chi$. Consequently, for any $[v]\in A(V_L^+)(\alpha)$, there exists $u\in M(1)^+$ such that $[v]=[u]*[B_{\alpha}]$ on $W$, and on $V_L^{T_{\chi},+}(0)$ for any $\chi$.
\end{lem} 
\begin{proof}By Theorem \ref{spangen}, we can conclude that there exist $u\in M(1)^+$, and $\beta\in L$ such that $\<\beta,\beta\><0$, $\<\alpha,\beta\><0$ and $$[B_{\alpha}]=[u]*[B_{\alpha+2\beta}].$$ Since $[u]$ acts as a constant on $M(1)(\theta)^+(0)$, says $p$, we have $[B_{\alpha}]=p[B_{\alpha+2\beta}]$ on $W$ and on $V_L^{T_{\chi},+}(0)$ for any $\chi$. On the other hand, 
$$[B_{\alpha}]=[B_{\alpha}]*[B_{2\beta}]=\epsilon(\alpha,2\beta)[B_{\alpha+2\beta}]=[B_{\alpha+2\beta}]$$ on $V_L^{T_{\chi},+}(0)$ for any $\chi$. Hence, $p=1$ and $[B_{\alpha}]=[B_{\alpha+2\beta}]$ on $W$.
\end{proof}
\begin{coro}\label{cb2action} Let $\alpha\in L$ such that $\<\alpha,\alpha\>=0$ and let $\rho\in L$ such that $\<\rho,\rho\><0$, $\<\alpha,\rho\><0$. Then $[B_{\alpha}]=[B_{\alpha+2\rho}]$ on $W$ and on $V_L^{T_{\chi},+}(0)$ for any $\chi$.
\end{coro} 
\begin{thm}\label{thm3} Let $W$ be an irreducible $A(V_L^+)$-modules such that $A^tW= 0$. Then there exists an irreducible $\hat{L}/K$-module $T_{\chi}$ with central character $\chi$ such that $W\cong T_{\chi}=V_L^{T_{\chi},+}(0)$.
\end{thm}
\begin{proof} First, we will show that for $\alpha,\beta\in L$, $[B_{\alpha}]*[B_{\beta}]=\epsilon(\alpha,\beta)[B_{\alpha+\beta}]$ on $W$. Let $\alpha,\beta\in L-\{0\}$.
\begin{description}
\item[case 1:] $\<\alpha,\beta\><0$. 

Then $[B_{\alpha}]*[B_{\beta}]\in A(V_L^+)(\alpha+\beta)$. In fact, $[B_{\alpha}]*[B_{\beta}]=[u]*[B_{\alpha+\beta}]$ for some $u\in M(1)^+$. By following the proof of Lemma \ref{baction} {\bf case 1}, we can show that $[B_{\alpha}]*[B_{\beta}]=\epsilon(\alpha,\beta)[B_{\alpha+\beta}]$ on $W$.

\item[case 2:] $\<\alpha,\beta\>>0$. Note that for this case $[B_{\alpha}]*[B_{\beta}]=[u]*[B_{-\alpha+\beta}]$ for some $u\in M(1)^+$.

\begin{description}
\item[case 2.1] $\<\alpha,\alpha\>> 0$, and $\<\beta,\beta\> \in 2\Z$. 

We have 
\begin{eqnarray*}
[B_{\alpha}]*[B_{\beta}]&=&\epsilon^{-1}(-\alpha,\alpha+\beta)[B_{\alpha}]*([B_{-\alpha}]*[B_{\alpha+\beta}])\\
&=&\epsilon^{-1}(-\alpha,\alpha+\beta)\epsilon(\alpha,-\alpha)[B_{\alpha+\beta}]\\
&=&\epsilon(\alpha,\beta)[B_{\alpha+\beta}],
\end{eqnarray*}
and 
\begin{eqnarray*}
[B_{\beta}]*[B_{\alpha}]&=&\epsilon^{-1}(\beta+\alpha,-\alpha)[B_{\alpha+\beta}]*[B_{-\alpha}]*[B_{\alpha}]\\
&=&\epsilon^{-1}(\beta+\alpha,-\alpha)\epsilon(-\alpha,\alpha)[B_{\alpha+\beta}]\\
&=&\epsilon(\beta,\alpha)[B_{\alpha+\beta}]
\end{eqnarray*}
on $W$.
\item[case 2.2:] $\<\alpha,\alpha\><0$, and $\<\beta,\beta\>< 0$.

By following the proof of Lemma \ref{baction} cases 2.1-2.3, we can show that $[B_{\alpha}]*[B_{\beta}]=\epsilon(\alpha,\beta)[B_{\alpha+\beta}]$ on $W$.
\item[case 2.3:] $\<\alpha,\alpha\><0$ and $\<\beta,\beta\>=0$. 

We have 
\begin{eqnarray*}
[B_{\alpha}]*[B_{\beta}]&=&[u]*[B_{-\alpha+\beta}]\\
&=&[u]*[B_{2\alpha}]*[B_{-\alpha+\beta}]\\
&=&[u]*\epsilon(2\alpha,-\alpha+\beta)[B_{\alpha+\beta}]\ \  \text{(by {\bf case 2.2})}\\
&=&[u]*[B_{\alpha+\beta}].
\end{eqnarray*}
Similarly, we have $[B_{\beta}]*[B_{\alpha}]=[u]*[B_{\alpha+\beta}]$. As in Lemma \ref{baction} {\bf case 1}, we have that $[B_{\alpha}]*[B_{\beta}]=\epsilon(\alpha,\beta)[B_{\alpha+\beta}]$ and $[B_{\beta}]*[B_{\alpha}]=\epsilon(\beta,\alpha)[B_{\alpha+\beta}]$ on $W$.
\item[case 2.4] $\<\alpha,\alpha\>=0$ and $\<\beta,\beta\>=0$.

We have
\begin{eqnarray*}
[B_{\alpha}]*[B_{\beta}]&=&[u]*[B_{-\alpha+\beta}]\\
&=&[u]*([B_{-\alpha+\beta}]*[B_{2\alpha-4\beta}])\\
&=&[u]*[B_{\alpha-3\beta}]\\
&=&[u]*[B_{\alpha-3\beta}]*[B_{-2\alpha+2\beta}]\\
&=&[u]*[B_{\alpha+\beta}].
\end{eqnarray*}
 By following the proof of Lemma \ref{baction} {\bf case 1}, we can show that $[B_{\alpha}]*[B_{\beta}]=\epsilon(\alpha,\beta)[B_{\alpha+\beta}]$ on $W$.
\end{description}
\item[case 3:] $\<\alpha,\beta\>=0$.

If one of the followings holds: $\<\alpha,\alpha\>> 0$, and $\<\beta,\beta\> \in 2\Z$, or $\<\alpha,\alpha\><0$, and $\<\beta,\beta\>< 0$ or $\<\alpha,\alpha\><0$ and $\<\beta,\beta\>=0$, then as in {\bf case 2}, we have that $[B_{\alpha}]*[B_{\beta}]=\epsilon(\alpha,\beta)[B_{\alpha+\beta}]$ and $[B_{\beta}]*[B_{\alpha}]=\epsilon(\beta,\alpha)[B_{\beta+\alpha}]$ on $W$.

Next, we consider the case when $\<\alpha,\alpha\>=0$ and $\<\beta,\beta\>=0$. Since there exists $\gamma_1\in L$ such that $\<\gamma_1,\gamma_1\><0$ and $\<-2\alpha+2\beta,\gamma_1\><0$, it implies that
\begin{eqnarray*}
[B_{\alpha+\beta}]&=&[B_{\alpha+\beta}]*[B_{-2\alpha+2\beta+2\gamma_1}]\\
&=&[B_{-\alpha+3\beta+2\gamma_1}]\\
&=&\epsilon^{-1}(-\alpha+3\beta,2\gamma_1)[B_{-\alpha+3\beta}]*[B_{2\gamma_1}]\\
&=&[B_{-\alpha+3\beta}].
\end{eqnarray*}
Similarly, since there exists $\gamma_2\in L$ such that $\<\gamma_2,\gamma_2\><0$ and $\<2\alpha-4\beta,\gamma_2\><0$, it follows that
\begin{eqnarray*}
[B_{\alpha}]*[B_{\beta}]&=&[u]*[B_{-\alpha+\beta}]\\
&=&[u]*([B_{-\alpha+\beta}]*[B_{2\alpha-4\beta+2\gamma_2}])\\
&=&[u]*[B_{\alpha-3\beta+2\gamma_2}]\\
&=&[u]*\epsilon^{-1}(\alpha-3\beta,2\gamma_2)[B_{\alpha-3\beta}]*[B_{2\gamma_2}]\\
&=&[u]*[B_{\alpha-3\beta}]\\
&=&[u]*[B_{\alpha+\beta}]
\end{eqnarray*}
on $W$ and on $V_L^{T_{\chi},+}(0)$ for any $\chi$. By following the proof of Lemma \ref{baction} {\bf case 1}, we can show that $[B_{\alpha}]*[B_{\beta}]=\epsilon(\alpha,\beta)[B_{\alpha+\beta}]$ on $W$.
\end{description}
Hence, for any $\alpha,\beta\in L$, we have $[B_{\alpha}]*[B_{\beta}]=\epsilon(\alpha,\beta)[B_{\alpha+\beta}]$ on $W$.

Next, to show that $W$ is isomorphic to some irreducible $\hat{L}/K$-module, we just need to follow the proof of Theorem \ref{cam1}.
\end{proof}

Next, we let $M$ be an irreducible $A(V_L^+)$-module such that $A^tM\neq 0$. 
\begin{lem} $M=A^tM$ and $M$ is a direct sum of $M(1)(\theta)^-(0)$.
\end{lem}
\begin{proof}  We first note that $M$ contains a simple $A^t$-modules which is isomorphic to $M(1)(\theta)^-(0)$ since $A^t$ is a simple algebra and $A^tM\neq 0$. By following the proof in Lemma \ref{am+act} and using the fact that $A(V_L^+)$ is spanned by $A(V_L^+)(\alpha)$ for all $\alpha\in L$ such that $\<\alpha,\alpha\>\neq 0$, we can show that $M^0=\{m\in M|A^tm=0\}$ is an $A(V_L^+)$-module. Moreover, we have that $M=A^tM$ and $M$ is a direct sum of $M(1)(\theta)^-(0)$.
\end{proof}

We will show that $M$ is isomorphic to $V_L^{T_{\chi,-}}(0)$ for some central character $\chi$.

For $\alpha\in L$ such that $\<\alpha,\alpha\>\neq 0$, we define
$$\tilde{B}_{\alpha}=2^{\<\alpha,\alpha\>-1}(E^{\alpha}-\frac{2\<\alpha,\alpha\>}{2\<\alpha,\alpha\>-1}E^t_{11}*E^{\alpha}),$$
where $E^t_{11}$ is defined with respect to an orthonormal basis $\{h_a|1\leq a\leq d\}$ of $\h$ such that $h_1\in \C\alpha$, and we set $$\tilde{B_0}=1.$$ Next, we let $\alpha\in L-\{0\}$ such that $\<\alpha,\alpha\>=0$. Then there exist $\gamma,\beta\in L_{\Q}$ such that $\alpha=\gamma+\beta$, $\<\gamma, \beta\>=0$, $\<\gamma, \gamma\>>0$ $,\<\beta,\beta\><0$. We define an action of $\tilde{B}_{\alpha}$ on $M$ in the following way: for $1\leq c\leq d$,
\begin{eqnarray*}
[\tilde{B}_{\alpha}]h_c(-1/2)&=&\frac{1}{2}[E^{\alpha}]h_c(-1/2)+\left(\frac{\<\gamma,\gamma\>}{1-2\<\gamma,\gamma\>}\delta_{1,c}-\frac{1}{2}\delta_{2,c}\right)[E^t_{11}]*[E^{\alpha}]h_c(-1/2)\\
& &+\left(\frac{\<\beta,\beta\>}{1-2\<\beta,\beta\>}\delta_{2,c}-\frac{1}{2}\delta_{1,c}\right)[E^t_{22}]*[E^{\alpha}]h_c(-1/2).\end{eqnarray*}
Here $E^t_{11}$ and $E^t_{22}$ are defined with respect to an orthonormal basis $\{h_a|1\leq a\leq d\}$ so that $h_1\in\C\gamma$, and $h_2\in \C\beta$, respectively.

Note that for $\alpha\in L$ such that $\<\alpha,\alpha\><0$, we have $[\tilde{B}_{2\alpha}]=1$ on $M$.
\begin{lem}\label{comm2} Let $\alpha\in L$ such that $\<\alpha,\alpha\>\neq 0$. Then for any $1\leq a,b\leq d$, $[\tilde{B}_{\alpha}]$ and $[E^t_{ab}]$ commute on $M$. Consequently, $[\tilde{B}_{\alpha}]$ commutes with the action of $[u]$ where $u$ is any member of $M(1)^+$.
\end{lem}
\begin{proof} The proof is similar to the proof of Lemma \ref{comm}.\end{proof}

\begin{lem}\label{lb3action} Let $\alpha\in L$ such that $\<\alpha,\alpha\>=0$. Then there exists $\beta\in L$ such that $\<\beta,\beta\><0$, $\<\alpha,\beta\><0$ and $[\tilde{B}_{\alpha}]=[\tilde{B}_{\alpha+2\beta}]$ on $M$, and $V_L^{T_{\chi},-}(0)$ for any $\chi$. Consequently, 
\begin{enumerate}
\item for any $[v]\in A(V_L^+)(\alpha)$, there exists $u\in M(1)^+$ such that $[v]=[u]*[\tilde{B}_{\alpha}]$ on $M$, and $V_L^{T_{\chi},-}(0)$ for any $\chi$. 
\item $[\tilde{B}_{\alpha}]$ commutes with the action of $[w]$ where $w$ is any member of $M(1)^+$. In particular, for any $1\leq a,b\leq d$, $[\tilde{B}_{\alpha}]$ and $[E^t_{ab}]$ commute on $M$. 
\end{enumerate}
\end{lem} 
\begin{proof} The proof is similar to Lemma \ref{lb2action}.\end{proof}
\begin{lem} For $\alpha,\beta\in L$, we have $[\tilde{B}_{\alpha}]*[\tilde{B}_{\beta}]=\epsilon(\alpha,\beta)[\tilde{B}_{\alpha+\beta}]$ on $M$.
\end{lem}
\begin{proof} By following the first part of the proof of Lemma \ref{modl}, we can show that when $\<\alpha,\beta\><0$, $[\tilde{B}_{\alpha}]*[\tilde{B}_{\beta}]=\epsilon(\alpha,\beta)[\tilde{B}_{\alpha+\beta}]$ on $M$. Next, by modifying the proof of Theorem \ref{thm3}, we will obtain that $[\tilde{B}_{\alpha}]*[\tilde{B}_{\beta}]=\epsilon(\alpha,\beta)[\tilde{B}_{\alpha+\beta}]$ on $M$ when $\<\alpha,\beta\>\geq 0$.
\end{proof}

\begin{thm}\label{thm4} Let $M$ be an irreducible $A(V_L^+)$-modules such that $A^tM\neq 0$. Then there exists an irreducible $\hat{L}/K$-module $T_{\chi}$ with central character $\chi$ such that $M\cong \h(-1/2)\otimes T_{\chi}=V_L^{T_{\chi},-}(0)$.
\end{thm}
\begin{proof}  The proof is the same as the proof of Theorem \ref{thm2}.\end{proof}
The following proposition is a consequence of Theorem \ref{thm3} and Theorem \ref{thm4}.
\begin{prop} Let $L$ be an even lattice of a finite rank equipped with a symmetric nondegenerate $\Q$-valued $\Z$-bilinear form that is neither positive or negative definite. Then any irreducible admissible $V_L^+$-module is isomorphic to $V_L^{T_{\chi},\pm}$ for any irreducible $\hat{L}/K$-module $T_{\chi}$ with central character $\chi$ such that $\chi(\iota(\kappa))=-1$.
\end{prop}

\end{document}